\documentstyle[12pt]{amsart}

\newcommand{\R}{{\mathbb R}}

\newcommand{\Z}{{\mathbb Z}}

\newcommand{\half}{{\frac{1}{2}}}

\newtheorem{theo}{{\sc Theorem}}[section]
\newtheorem{defn}[theo]{{\sc Definition}}
\newtheorem{cor}[theo]{{\sc Corollary}}

\newtheorem{lem}[theo]{{\sc Lemma}}
\newtheorem{prop}[theo]{{\sc Proposition}}
\newenvironment{proof}{\noindent{\em Proof:\/}}{\qed \medskip}

\title[Norms of modes and quasi-modes revisited ]{Norms of modes and quasi-modes revisited}

\author{John A. Toth and Steve Zelditch}
\address{Department of Mathematics and Statistics, McGill University,
Montreal, CANADA, H3A-2K6}
\address{Department of Mathematics, Johns Hopkins University, Baltimore, MD
21218, USA}

\thanks{\\ Research partially supported by NSERC grant \#OGP0170280\\ Research
partially supported by  NSF grant \#DMS-0100474}

\date{November 19, 2001}

\begin{document}

\maketitle

\addtolength{\baselineskip}{1pt}

\begin{abstract}  This article is devoted to the analysis of
eigenfunctions (modes) and approximate eigenfunctions
(quasi-modes) of the Laplacian on a compact manifold $(M, g)$ with
completely integrable geodesic flow. We give a new proof of the
main result of \cite{TZ} that $(M, g)$ with integrable Laplacians
and  with uniformly bounded eigenfunctions must be flat. The proof
is based on the use of Birkhoff normal forms and on  a comparison
of modes and quasimodes. In the process,
 we discuss   tunnelling between resonant tori and give a proof that
eigenfunctions concentrate on individual tori in the non-resonant
case.   We also give brief expositions of results in \cite{TZ2,
SZ}.

\end{abstract}

\setcounter{page}{1} \setcounter{section}{-1}

\section{Introduction}  In this article we are interested in the
semiclassical $h \to 0$   asymptotics of modes (eigenfunctions)
$\phi_j(h)$
\begin{equation} \label{EVP} (h^2 \Delta + V) \phi_j = \lambda_j^2(h)  \phi_j, \;\;\;
\langle \phi_j, \phi_k \rangle = \delta_{jk} \end{equation} of
Schroedinger operators on  a Riemannian manifold $(M, g)$. Usually
we take $(M, g)$ to be compact for convenience, but equally
interesting problems occur on non-compact manifolds such as
$\R^n$.  Many questions about eigenfunctions are motivated by
modern mesoscopic physics and quantum chaos, in particular by the
physics of quantum dots, wires and corrals (for a recent sample,
see e.g. \cite{FE,FHHMLE,K,  Mi, Mi2}).  Let us mention just a few
topics and  recent articles from  the mathematics and physics
literature :
\begin{itemize}

\item Norms of eigenfunctions: $||\phi_j||_{L^p}$ \cite{ ABST, AS, Bour, D,  IS, SZ};

\item Eigenfunction scarring \cite{AA, Be2, K,  RS, W};

\item Volume  and distribution of nodal surfaces $\{\phi_j = 0\}$ \cite{BGS, BS, CHENG, DF, DF2, JL};

\item Number and distribution of critical points $\{\nabla \phi_j = 0\}$ \cite{JN, SZ};

\item Matrix elements $\langle A \phi_j, \phi_k \rangle$ of
pseudodifferential operators \cite{KH, Z};

\item Resonance and tunnelling between wells or  tori \cite{HS, HSu, M, Wi}.

\end{itemize}

Most of the classical and more recent  problems about
eigenfunctions are believed to be very difficult.
 The main reason for the difficulty is that in
general there are no explicit formulae for eigenfunctions or even
for approximations to them. Lacking such explicit formulae, one
has to rely on local approximations that lose track of the global
geometry of $(M, g)$ or on  wave equation or resolvent methods, which
involve averages over many eigenfunctions. It is thus difficult to
obtain global information about individual modes.

In this article, we would like to review and revisit some recent
results   about $L^p$-norms of eigenfunctions, obtained in part by
the authors and in part by one of the authors  in collaboration
with C. Sogge. For some $(M, g)$, $L^p$- norms of eigenfunctions
are extremally large \cite{SZ} and for some they are extremally
small \cite{TZ, TZ2}. The problems studied in \cite{TZ, TZ2, SZ}
were to characterise metrics and potentials producing extremal
growth rates of eigenfunctions. Of course,  one would like to have
some control over growth rates of eigenfunctions for any $(M, g,
V)$. Our main focus in this article  is to give a different and in
some ways more natural proof of the
 results on $L^{\infty}$ norms of modes in the quantum integrable case
  of \cite{TZ}, based on the construction of  quasimodes.

   On a
  more expository level, we would also like to advertise
 quantum integrable systems as a simple setting for harmonic
 analysis. Quantum integrable Laplacians are generalizations of
Laplacians on  flat tori and of  harmonic oscillators, and have
most of the features which make  these settings attractive for
harmonic analysis. Yet they also bring in interesting new
geometric and analytic phenomena.  Most
 of the special eigenfunctions studied  in the  literature are quantum integrable, see
 e.g. \cite{B.K.S, Bour, CV1, D, JN, KMS}.  To make the ideas more
 accessible to  harmonic analysts, we describe a number of
 two-dimensional examples in some detail in Section
 (\ref{EXAMPLES}).

To motivate the problem and result, let us consider the setting
where we first learn about eigenfunctions, namely  the flat torus
$\R^n/ L$ where $L$ is lattice such as $\Z^n.$ The eigenfunctions
$\phi_{\lambda}$ are of course the exponentials $e^{i \langle
\lambda, x \rangle}$ where $\lambda \in L^*$, the dual lattice.
The harmless looking facts
$$\phi_{\lambda}(x) = e^{i \langle \lambda, x \rangle},\;\;\;
|\phi_{\lambda}(x)| = 1, \;\; (\forall x)$$ about the most
familiar eigenfunctions raise the questions: when can
eigenfunctions be expressed in this kind of WKB form $a(x) e^{i
S(x)/h}$ or as an oscillatory integral?  And how common is it that
eigenfunctions have uniformly bounded sup norms?

Heuristic  answers to these and other questions about
eigenfunctions are suggested by correspondence principles between
excited states of quantum systems and their underlying classical
dynamics. The quantum Hamiltonian is simply the Laplacian
$\Delta$, while its classical counterpart is its bicharacteristic
flow, namely the geodesic flow $G^t: S^*M \to S^*M.$ The clearest
statements concern the most extreme dynamical regime, integrable
systems and chaotic systems, with mixed and somewhat fuzzy
statements about intermediate systems. The following description
summarizes the current conjectural picture of eigenfunctions; it
is not intended to be exhaustive. Eigenfunctions are described
according to the regime of the corresonding classical dynamics.

\begin{itemize}

\item The Integrable regime:  Approximate  eigenfunctions may be constructed by a
WKB method associated to invariant tori satisfying Bohr-Sommerfeld
quantization conditions. Sequences of  eigenfunctions localize in
phase space on (unions of ) the invariant tori for the geodesic flow. They `scar'
(concentrate singularly) on the projections of those tori.

\item The intermediate KAM regime: WKB quasimodes may be constructed for a positive measure
set of ``Cantori". There are regular and irregular modes corresponding to the
Cantori and the regions of instability between them.

\item The chaotic regime: eigenfunctions behave like random waves, e.g. like
random spherical harmonics. They do not localize in phase space, but rather become
uniformly distributed, as do their critical points and zero sets. Only low-lying
modes should scar.

\end{itemize}

  Quasi-modes are
essentially parametrices for modes, in some ways analogous to
parametrices for operator kernels (e.g. the Hadamard parametrix
for the heat or wave kernel). Unlike operator parametrices,
however, they do not necessarily approximate individual modes. We
return to this point below.

Quasi-modes are usually thought of as oscillatory integrals
associated to geodesic-flow invariant Lagrangean (or isotropic)
submanifolds of an energy level $\{H(x, \xi) = |\xi|^2_g + V(x) =
E\}  \subset T^*M$.  There are axioms for quasimodes \cite{CV2}
which allow for more general constructions,  in particular in the
KAM case \cite{CV2, L, P}. Existence of quasimodes associated to a
Lagrangean submanifold $\Lambda$ involves a Bohr-Sommerfeld
quantization condition on $\Lambda$ and some {\it stability}
properties. For instance, there are always quasi-modes associated
to stable closed geodesics \cite{CV2, R}, but there are no known
quasi-mode constructions associated to unstable closed geodesics.
For this reason, quasi-mode constructions are most effective in
integrable or KAM systems, and these regimes  represent the best
chance of understanding individual eigenfunctions. In the opposite
chaotic regime, there are no known quasi-modes and heuristic
reasons suggest there aren't any explicitly constructible
approximate solutions. Indeed, the modes are believed to behave
like random waves and to be rather featureless. Finally, in the no
man's land between the nearly integrable and the quite chaotic,
there are few if any heuristics or results of a global nature. The
only known results are obtained by piecing together local
information on eigenfunctions.

Let us now state some recent rigorous results. First we need some
notation: by $\Psi^m(M)$ we denote  the space of mth order
pseudodifferential operators over $M$. We say that an operator $P
\in \Psi^1(M)$ is quantum integrable if there exist $P_2, \dots,
P_m \in \Psi^1(M)$ such that $[P_j, P_k] = 0$ and such that $P_1^2 + \cdots +
P_m^2$ is elliptic. We denote \ by $p_k = \sigma_{P_k}$
 the principal symbol of $P_k$.  Since $\sigma_{[P_i, P_j]} =
\{p_i, p_j\}$ (the Poisson bracket), it follows that the $p_j$'s
generate a homogeneous
 Hamiltonian action $\Phi_t$ of $t \in \R^n$ on $T^*M - 0$ with moment map
$${\cal P} : T^*M - 0 \rightarrow \R^n,\;\;\;\;\;\;{\cal P} = (p_1, \dots,
p_n).$$ We denote the image ${\cal P}(T^*M - 0)$ by $B$, by
$B_{reg}$ (resp. $B_{sing}$)  the regular values (resp. singular
values)  of the moment map.

We are going to impose a `finite complexity' assumption on the
classical integrable system. Suffice it to say that it holds for
all systems of interest in physics:  for each
$b=(b^{(1)},...,b^{(n)})\in B$, let $m_{cl}(b)$ denote the number
of ${\Bbb R}^{n}$-orbits of the joint flow $\Phi_{t}$ on the level
set ${\cal P}^{-1}(b)$. Then
\begin{equation} \label{HYP} \mbox{\bf Finite complexity condition}:\;\;\;  \exists M: \;  m_{cl}(b) < M\;\;\;(\forall b \in B).  \;\;  \end{equation}
When $b \in B_{reg}$, then ${\cal P}^{-1}(b)$ is the union of
$m_{cl}(b)$ isolated  Lagrangian tori. If $b \in B_{sing}$, then
${\cal P}^{-1}(b)$  consists of a finite number of connected
components, each of which is a finite union of orbits. These
orbits may be Lagrangian tori, singular compact tori (ie. compact
tori of dimension $<n$) , or non-compact orbits consisting of
cylinders or planes.

On the quantum level we will assume:
\begin{equation} \label{QHYP} \begin{array}{lllll} \mbox{\bf Bounded eigenvalue multiplicity }: & \exists M': \;m(\lambda) \leq M'\;\;  & &
 (\forall \lambda; \;m(\lambda) = dim V_{\lambda}). \end{array} \end{equation}
This reduces estimates of eigenfunctions to estimates of joint eigenfunctions
of the commuting operators. When unbounded  multiplicities occur, such as on $\R^n/\Z^n$,
then it is well-known that one may construct sequences of eigenfunctions with unbounded sup norms even on flat tori.

\begin{theo} \label{RM}  Suppose that $\Delta$ is a quantum completely
integrable Laplacian on a compact Riemannian manifold $(M, g)$,
and suppose that the  corresponding moment map satisfies
(\ref{HYP}). Then: if the $\Delta$-eigenfunctions have uniformly
bounded $L^{\infty}$ norms, then $(M, g)$ is flat.

More generally, suppose that $\hbar^2 \Delta + V$ is a quantum
completely integrable Schroedinger operator,  and that the
corresponding moment map ${\cal P}$  is proper and satisfies
(\ref{HYP}).   Assume there exists an energy level $E$ such that
the eigenfunctions have uniformly bounded $L^{\infty}$ norms as $h
\to 0$. Then: if  $E > \max V,$ and $(M,(E - V) g)$ is flat. If
(a) (or (b)) holds for all energy levels $E$ in an interval $E_1 <
E < E_2$, then $(M, g)$ is flat and $V$ is constant.

\end{theo}

To prove the result, we must  relate the analysis of modes to the
geometry of the Liouville foliation (by invariant tori) of $T^*M$
of the completely integrable system. In \cite{TZ}, this relation
was brought out through trace formulae. In this article, we
approach it  through quasimodes. With some work, we believe the
results should generalize to some degree to quantum KAM systems,
i.e. to small perturbations of completely integrable systems.
However, we leave that extension for the future.

Roughly speaking, the idea is this: Quasi-modes are directly
associated to the individual Lagrangian tori of the completely
integrable system. Hence we can directly relate norms of
quasimodes to geometry of the tori.  To draw conclusions about
norms of modes, we must then relate modes and quasimodes.  Since
the classic article
 ``Modes and Quasi-modes"  \cite{A} of V. I. Arnol'd (also  based
on  examples due to  A. I. Schnirelman), attention has been drawn
to the danger of of drawing conclusions about modes from
conclusions about quasimodes. As Arnold cautioned, modes are not
always approximated by quasimodes, but rather may be linear
combinations of quasimodes with close quasi-eigenvalues. If the
linear combinations involve `many' closely spaced
quasi-eigenvalues, then the connection between modes and
quasi-modes becomes blurred. However, this danger mainly
pertains to scalar spectra, i.e. the spectrum of one operator;
it does not occur in the approximation of joint eigenfunctions of
a completely integrable set of operators by joint quasimodes,
 since the joint spectrum
in $\R^n$ is much better separated than the scalar spectrum of the
Hamiltonian alone. For the joint spectrum, a simple non-resonance
condition is sufficient to ensure that modes are even
well-approximated by individual quasi-modes.

When modes of integrable systems  can be connected with invariant
tori for the geodesic flow, one can hope to relate growth rates of
$L^p$-norms of sequences of eigenfunctions to singularities of
projections of the associated tori. For instance, as will be
discussed further in Section (\ref{EXAMPLES}), invariant
eigenfunctions on surfaces of revolution  correspond to to
invariant torus in $S^*S^2$ consisting of meridian geodesics
between the poles. This torus has a ``blow-down'' singularity
under the natural projection $\pi: S^*S^2 \to S^2$, i.e. the fiber
$S^*_{NP}S^2 \to \{NP\}$ over the north pole is blown-down to a
point under the projection. This kind of singularity causes the
largest concentration of the mass of a mode at a point; hence the
invariant eigenfunctions which correspond to the meridian torus
are (at least in order of magnitude) extremals for the
$L^{\infty}$ and other  $L^p$ norms for sufficiently large $p$.
Other important  tori in the $L^p$ theory are the {\it singular
tori} where the torus drops dimension, e.g. the equators of a
surface of revolution. Eigenfunctions concentrated on such
singular tori are (in order of magnitude) extremals for low $L^p$
norms: although they do not blow up as much, they blow up on a set
of larger measure and hence are favorable for $L^p$ norm blow-up.

Let us now recall some other recent results; we refer to \cite{TZ2, SZ} for
further discussion.

The next result also pertains to quantum completely integrable systems. We impose
a stronger hypothesis on the classical system, namely Eliasson non-degeneracy
(see \cite{TZ2} for background). We then obtain quantitative results on the $L^p$
 blow-up of eigenfunctions in non-flat cases, i.e. in cases where (by the previous
theorem) some blow-up must occur. We emphasize that our notation
for eigenvalues (\ref{EVP})  is such that $\Delta \phi_k =
\lambda_k^2 \phi_k$.

\begin{theo} \label{LP}  Suppose that $(M, g)$ is a compact  Riemannian manifold
with completely integrable geodesic flow  satisfying Eliasson's
non-degeneracy condition. Then $(M, g)$ is non-flat if and only if
for every $\epsilon >0$,  there exists a sequence of
$L^2$-normalized eigenfunctions $\phi_k$ satisfying:
$$\left\{ \begin{array}{l}   \| \phi_{k} \|_{\infty} \geq C(\epsilon) \lambda_{k}^{\frac{1}{4} -
\epsilon}. \\ \\
 \| \phi_{k} \|_{p} \geq C(\epsilon) \lambda_{k}^{ \frac{ p - 2}{ 4p } - \epsilon }\end{array} \right.$$
   \end{theo}

Modulo the $\epsilon$'s, these estimates are sharp. They are based
on the existence of codimension one singular leaves of the
Lagrangian foliation, and on the construction of precise
quasimodes attached to non-degenerate singular leaves. In general,
a completely integrable system does not have higher codimension
singular leaves (e.g. one may take the product $M = {\bf T}^{n-2}
\times S^2$ of a flat torus and a surface of revolution to obtain
completely integrable  $(M, g)$ of any dimension whose only
singular leaves are of codimension one). When higher codimension
leaves occur, the blow-up rate increases. In the case where $(M,
g)$ has a one-dimensional singular leaf , then our methods show
that there exists  a sequence of eigenfunctions satisfying: $   \|
\phi_{k} \|_{\infty} \geq C(\epsilon) \lambda_{k}^{\frac{n-1}{4} -
\epsilon}. $ Such a one-dimensional leaf is simply a stable
elliptic geodesic $\gamma$ which is fixed by the other Hamiltonian
flows of the completely integrable system. This estimate is
comparable to the recent estimate of Donnelly \cite{D} on the
blow-up rate of a sequence of eigenfunctions on certain $(M, g)$
with an isometric $S^1$ action. Note that Donnelly denotes the
eigenvalue by $\lambda$ while we denote it by $\lambda^2$.
Donnelly also assumes the existence of a non-degenerate singular
orbit $\gamma$ of the joint action of the geodesic flow and $S^1$
action.  Such a $\gamma$ is a a stable closed geodesic, such as
the equator on a convex surface of revolution,  which is invariant
under the $S^1$ action. The corresponding sequence of
eigenfunctions is actually a sequence of singular quasi-modes
associated to $\gamma$. Thus, it is the same phase space
phenomenon in either setting which is responsible for this common
blow-up rate.

We now turn to results by one of the authors with C. Sogge on the opposite extreme:
which $(M, g)$ have maximal eigenfunction growth?  It follows from classical remainder estimates
of Avakumovic and Hormander that $L^2$-normalized eigenfunctions
$\Delta \phi_{\lambda} = \lambda^2 \phi_{\lambda}$ have $L^{\infty}$-norms satisfying
$$||\phi_{\lambda}||_{\infty}  \leq C \lambda^{\frac{n-1}{2}},\;\;\; n = \dim M.$$
When a sequence of eigenfunctions of $(M, g)$ actually achieves this bound, we say
that $(M, g)$ has {\it maximal eigenfunction growth}.
 The  main result of \cite{SZ}   implies  a necessary condition on a compact Riemannian manifolds $(M, g)$ with maximal eigenfunction growth: there  must exist
a point $x \in M$ for which the set
\begin{equation}  {\mathcal L}_x = \{ \xi \in S^*_xM : \exists T: \exp_x T \xi = x\} \end{equation}
of directions of geodesic loops at $x$ has positive surface
measure.  Here, $\exp$ is the exponential map, and the measure
$|\Omega|$ of a set $\Omega$ is the one  induced by the metric
$g_x$ on $T^*_xM$. For instance, the poles $x_N, x_S$ of a surface
of revolution $(S^2, g)$ satisfy $|{\mathcal L}_x| = 2 \pi$. More generally,
we may consider maximal $L^p$ growth of eigenfunctions. The universal esimates
in this case are due to Sogge \cite{So1, So2}.

\begin{theo}\label{eigcor}
If  $|{\cal L}_x|=0$ for every $x\in M$ then
\begin{equation}\label{eig2}
\frac{\|\phi\|_{L^p(M)}}{\|\phi\|_{L^2(M)}} = o(\lambda^{\delta(p)}), \, \, \, \,
p>\tfrac{2(n+1)}{n-1}\, \end{equation} where
\begin{equation} \delta(p)=
\begin{cases}
n(\tfrac12-\tfrac1p)-\tfrac12, \quad \tfrac{2(n+1)}{n-1}\le p\le
\infty
\\
\tfrac{n-1}2(\tfrac12-\tfrac1p),\quad 2\le p\le
\tfrac{2(n+1)}{n-1}.
\end{cases}
\end{equation}

\end{theo}

It is  proved that  generic metrics $g$ on any $M$ satisfy the hypothesis of
having zero-measure loop directions at every point. Hence, maximal eigenfunction
growth is a rare condition.
The strongest conclusion is available for real analytic Riemannian manifolds:

\begin{theo} \label{ANAL0} Suppose that $(M,g)$ is a real analytic
with maximal eigenfunction growth.
Then there exists a point $m \in M$  such that all geodesics    issuing
from the point $m$ return to $m$ at the same time  $\ell$. In particular,
if dim $M$ = 2, then $M$ is topologically a 2-sphere $S^2$ or a
real projective plane $\R P^2$.  \end{theo}

It would be interesting to investigate maximal eigenfunction growth in small
$L^p$ norms. It is explained in \cite{So1, So2}  that the geometry of eigenfunctions changes at the
critical $p = \tfrac{2(n+1)}{n-1},$ and that sizes of loop directions is no longer
relevant. In the case of the found sphere, at  this value of $p$ the extremals for the $L^p$ norm (at least in
order of magnitude) shift from eigenfunctions which localize at a point (zonal spherical harmonics)  to eigenfunctions
which localize on a closed geodesic (highest weight spherical harmonics).
It is possible that maximal $L^p$ growth in the small range
of $p$ is related to existence of the same stable elliptic geodesics $\gamma$ which
gave rise to the $L^p$ estimates in \cite{TZ2, D}.

\section{Quantum integrable systems}

We begin by recalling the definition of a completely integrable system.

\subsection{Completely integrable systems}

By a completely integrable system on $T^{*}M$ we mean  a set of
$n$ independent,   $C^{\infty}$ functions  $p_1, \dots, p_n$,
 on $T^{*}M$ satisfying:

\medskip

\begin{tabular}{l} $\bullet \,\, \{ p_{i}, p_{j} \}=0$ \,\, for all $1 \leq i,j \leq n$;\\
$\bullet \,\, dp_{1} \wedge dp_{2} \wedge \cdot \cdot \cdot \wedge dp_{n} \neq 0$ \,\, on an open dense subset of $T^{*}M.$
\end{tabular}
\medskip

The associated moment map is defined by
\begin{equation} \label{MM} {\cal P} = (p_1, \dots, p_n): T^*M \rightarrow
B \subset\R^n. \end{equation}  We  refer to to the set $B$ as the `image of the moment map.'
The Hamiltonians generate an action of $\R^n$ defined by
$$ \Phi_t = \exp t_1 \Xi_{p_1} \circ \exp t_2 \Xi_{p_2} \dots \circ \exp t_n
\Xi_{p_n}.$$
We often denote  $\Phi_t$-orbits by   $\R^n \cdot (x, \xi) $. The isotropy group of $(x, \xi)$ will be denoted by
 ${\mathcal I}_{(x, \xi)}.$ When  $\R^n \cdot (x, \xi) $ is a compact Lagrangean orbit, then ${\mathcal I}_{(x, \xi)}$
is a lattice of full rank in $\R^n$, and is  known
as the `period lattice', since it consists of the `times' $T \in \R^n$ such that $\Phi_T |_{\Lambda^{(j)}(b)} = Id.$

We will need the following:
\begin{defn}\label{SING}
We say that:

\begin{itemize}

\item  $b \in B_{sing}$ if ${\cal P}^{-1}(b)$ is a singular level of the moment map, i.e. if
there exists a point $(x, \xi) \in {\cal P}^{-1}(b)$  with  $dp_{1} \wedge \cdot \cdot \cdot \wedge dp_{n}(x,\xi) = 0$.
Such a point $(x, \xi)$ is called a singular point of ${\cal P}$.

\item a connected component of ${\cal P}^{-1}(b)$ ($b \in B_{sing}$)is a singular component if it  contains a singular point ;

\item an orbit $\R^n \cdot (x, \xi) $ of $\Phi_{t}$ is  singular if it is non-Lagrangean, i.e.   has dimension $<n$;

\item $b \in B_{reg}$ and that ${\cal P}^{-1}(b)$ is a regular level if all points $(x, \xi) \in {\cal P}^{-1}(b)$ are regular, i.e. if $dp_{1} \wedge \cdot \cdot \cdot \wedge dp_{n}(x, \xi) \not= 0$.

\item   a component of  ${\cal P}^{-1}(b)$ ( $b \in B_{sing} \cup B_{reg}$) is regular if it contains no singular points.

\end{itemize}

\end{defn}

By the Liouville-Arnold theorem [AM], the orbits of the joint flow $\Phi_{t}$ are
diffeomorphic to $\R^k \times T^m$ for some $(k,m), k + m \leq n.$
By the properness assumption on  ${\cal P}$, a regular level  has the form
\begin{equation} \label{CI1}
{\cal P}^{-1}(b) = \Lambda^{(1)}(b) \cup \cdot \cdot \cdot \cup
 \Lambda^{(m_{cl})}(b) , \;\;\;(b \in B_{reg})
\end{equation}
\noindent where each $\Lambda^{(l)}(b) \simeq T^n$ is an $n$-dimensional Lagrangian
torus.  The  classical (or geometric) multiplicity function   $m_{cl} (b) = \# {\cal P}^{-1} (b)$,
i.e. the number
of  orbits on the level set ${\cal P}^{-1}(b)$, is constant on    connected components of $B_{reg}$ and the moment map
(\ref{MM}) is a fibration over each component with fiber (\ref{CI1}).
 In  sufficiently small neighbourhoods $\Omega^{(l)}(b)$ of each component torus,
$\Lambda^{(l)}(b)$, the Liouville-Arnold theorem also gives the existence of local action-angle variables
$(I^{(l)}_{1},...,I^{(l)}_{n}, \theta^{(l)}_{1},...,\theta^{(l)}_{n})$ in terms
of which the joint flow of $\Xi_{p_{1}},...,\Xi_{p_{n}}$ is linearized [AM]. For convenience, we henceforth normalize the action variables $I^{(l)}_{1},...,I^{(l)}_{n}$ so that $I^{(l)}_{j} = 0; \, j=1,...,n$ on the torus $\Lambda^{(l)}(b)$.

Quantum integrability involves  semiclassical pseudo-differential
operators, so let us pause to set up the notation.  Given an open
$U \subset {\Bbb R}^{n}$, we say that $a(x,\xi;\hbar) \in
C^{\infty}(U \times {\Bbb R}^{n})$ is in the symbol class
$S^{m,k}(U \times {\Bbb R}^{n})$, provided
$$ |\partial_{x}^{\alpha} \partial_{\xi}^{\beta} a(x,\xi;\hbar)| \leq
C_{\alpha \beta} \hbar^{-m} (1+|\xi|)^{k-|\beta|}.$$ \noindent We
say that $ a \in S^{m,k}_{cl}(U \times {\Bbb R}^{n})$ provided
there exists an asymptotic expansion:
$$ a(x,\xi;\hbar) \sim \hbar^{-m} \sum_{j=0}^{\infty} a_{j}(x,\xi)
\hbar^{j},$$ \noindent valid for $|\xi| \geq \frac{1}{C} >0$  with
$a_{j}(x,\xi) \in S^{0,k-j}(U \times {\Bbb R}^{n})$ on this set.
We denote the associated $\hbar$ Kohn-Nirenberg quantization by
$Op_{\hbar}(a)$, where this operator has Schwartz kernel given
locally by the formula:
$$Op_{\hbar}(a)(x,y) = (2\pi \hbar)^{-n} \int_{{\Bbb R}^{n}}
e^{i(x-y)\xi/\hbar} \,a(x,\xi;\hbar) \,d\xi.$$ By using a
partition of unity, one constructs a corresponding class,
$Op_{\hbar}(S^{m,k})$, of properly-supported
$\hbar$-pseudodifferential operators acting on $C^{\infty}(M)$.

\begin{defn} \label{SCQCI}
We say that the operators $Q_{j} \in Op_{\hbar}(S^{m,k}_{cl});
\,\,j=1,...,n$, generate a semiclassical quantum completely
integrable system if
$$ [Q_{i}, Q_{j}]=0; \,\,\, \forall{1 \leq i,j \leq n},$$
\noindent and the respective semiclassical principal symbols
$q_{1},...,q_{n}$ generate a classical integrable system with
$dq_{1} \wedge dq_{2} \wedge \cdot \cdot \cdot \wedge dq_{n} \neq
0$ on a dense open set $\Omega \subset T^{*}M-0$. We also assume
that the finiteness condition (\ref{HYP}) is satisfied.

\end{defn}

Let us pause to relate this notion to the possibly more familiar
notion of $n$ commuting (ordinary) pseudo-differential operators
$P_1, \dots, P_n$. Such operators have a joint spectrum. In
eigenvalue asymptotics for a joint spectrum, it is natural to
consider  ladders of eigenvalues:
\begin{equation} \label{QCI1}
\{ (\lambda_{1k},...,\lambda_{nk}) \in Spec(P_{1},...,P_{n});
\,\forall j=1,..,n, \,  \lim_{k\rightarrow
\infty}\frac{\lambda_{jk}}{|\lambda_{k}|}  \, = \, b_{j}\},
\end{equation}
\noindent where $|\lambda_{k}|:= \sqrt{ \lambda_{1k}^{2} + ... +
\lambda_{nk}^{2} }$.  We then introduce the semiclassical
parameter $\hbar$ with values in a sequence $\{ \hbar_{k}^{-1};
k=1,2,3,... \}$ with $\hbar_{k} \rightarrow 0$ as $k \rightarrow
\infty$. Now, consider the semiclassically scaled operators:
\begin{equation} \label{QCI2}
 Q_{j} := \hbar P_{j}; \,\,\,\,j=1,2,...,n.
\end{equation}
\noindent In the semiclassical regime, it is convenient to work
with the $Q_{j}$'s rather than the $P_{j}$'s.  For the
eigenfunctions $\phi_{b,\hbar}$ associated with the ladder
(\ref{QCI1}) we clearly have:
\begin{equation} \label{QCI3}
Q_{j} \phi_{b,\hbar} =  b_{j} \phi_{b,\hbar} + \mu_{jk}
\phi_{b,\hbar}, \,\,\,\mbox{where}\,\,\, \mu_{jk} = o(1)
\,\,\,\mbox{as} \,\hbar \rightarrow 0.
\end{equation}

\noindent The operators $Q_{j}; j=1,...,n$ are semiclassical
pseudodifferential operators.

Consequently, in the homogeneous case, the operators $P_{1} =
\sqrt{\Delta},P_{2},...,P_{n}$ are classical pseudodifferential
operators of order one, it readily follows that $Q_{j}:= \hbar
P_{j} \in Op(S^{0,1}_{cl})$ generate a semiclassical quantum
integrable system in the sense of Definition \ref{SCQCI}.

\section{Quasimodes and Birkhoff normal forms }

 First, some notation. Following the convention in \cite{CP},
we denote $\hbar$-microlocal equivalence  on an open set, $\Omega
\subset T^{*}M$ by $=_{\Omega}$. Also, recall that \cite{CP} a
family of distributions, $u_{\hbar}$, where $\hbar \in
(0,\hbar_{0}]$ is called admissible if for any cutoff function,
$\chi(x,\xi) \in C^{\infty}_{0}(\Omega)$, there is an $N\in {\Bbb
Z}$ and $s \in {\Bbb Z}$, such that
$$\| \hbar^{N} Op_{\hbar}(\chi) u_{\hbar} \|_{H^{-s}} = {\cal O}(1)$$
\noindent  uniformly for $\hbar \in (0,\hbar_{0}]$.

\subsection{Inverse Birkhoff normal form}

All of our results on localization of modes and quasimodes  are
based on the use of Birkhoff normal forms.
 We  begin by defining such normal forms.

 Let $b \in B_{reg}$ and consider the  regular level set ${\cal P}^{-1}(b)$ of the moment map  with
${\cal L}^{(l)}$ denote the pull-back of the Maslov line bundle
over the affine torus, ${\Bbb T}^{n}$, given by $I_{1}^{(l)}=
...=I_{n}^{(l)}=0$ and  parametrized by the angle variables
$\theta_{1}^{(l)},...,\theta_{n}^{(l)}$. We denote the space of
smooth sections of this bundle by $C^{\infty}({\Bbb T}^{n}; {\cal
L}^{(l)})$ . The image of $\Omega^{(l)}(b)$ under the map to
normal form will be denoted by $ \Omega_{0}^{(l)}= {\Bbb T}^{n}
\times D_{1}$, where $D_{1}$ denotes a sufficiently small
neighbourhood of $0 \in {\Bbb R}^{n}$.  We start with a quantum
normal-form construction which can be found in a somewhat
different form in \cite{CV2}.
\begin{lem}
 \label{QBNF}
 For  $l=1,...,k$ and $j=1,...,n$, there exist
$\hbar$-Fourier integral operators, $U^{(l)}_{b}: C^{\infty}(M)
\rightarrow C^{\infty}({\Bbb T}^{n}; {\cal L}^{(l)})$,
microlocally elliptic
 on $\Omega^{(l)}$, together with  $C^{\infty}$ symbols,
$f_{j}^{(l)}(x;\hbar) \sim \sum_{k=0}^{\infty} f^{(l)}_{jk}(x)
\hbar^{k}$, with $f_{j0}(0)=0$ such that:
$$U^{(l)*} f_{j}^{(l)}(Q_{1}-b_{1},...,Q_{n}-b_{n};\hbar)U^{(l)} =_{\Omega_{0}^{(l)}}
\frac{\hbar}{i} \frac{\partial}{\partial \theta_{j}}.$$ \noindent
Moreover, when $Q_{1},...,Q_{n}$ are formally self-adjoint, the
operator $U^{(l)}_{b}$ can be taken to be microlocally unitary.
\end{lem}
\noindent{\em Proof}: To simplify the writing a little, we drop
the superscript $(l)$ in the following and denote the $\hbar$
principal symbols of $Q_{1},...,Q_{n}$ by $q_{1},...,q_{n}$. Let
$\kappa: \Omega(b) \rightarrow {\Bbb T}^{n} \times D_{1}$ be a
canonical map with the property that,
$$ q_{j} \circ \kappa (\theta,I) = \tilde{q}_{j}(I),\,\,\,\,j=1,...,n,$$
\noindent for some $\tilde{q}_{j} \in C^{\infty}(D_{1})$. Then,
since $dq_{1},...,dq_{n}$ are linearly-independent on
$\Omega(b)$, by the inverse function theorem, there exists $f_{j0}
\in C^{\infty}({\Bbb R}^{n})$ with $f_{j0}(0)=0$  such that:
\begin{equation} \label{sec43}
f_{j0}( q_{1}-b_{1},...,q_{n}-b_{n}) \circ \kappa =
I_{j},\,\,\,\,j=1,...,n.
\end{equation}
\noindent By the semiclassical Egorov Theorem, there exists an
$\hbar$-Fourier integral operator $U_{0}: L^{2}({\Bbb T}^{n};{\cal
L}) \rightarrow L^{2}(M)$ quantizing $\kappa$ such that for
$j=1,...,n$,
\begin{equation} \label{sec44}
U_{0}^{-1} f_{j0}(Q_{1}-b_{1},...,Q_{n}-b_{n}) U_{0} =_{{\Bbb
T}^{n} \times D_{1}} Op_{\hbar}(I_{j}) + \hbar
\,Op_{\hbar}(r_{1}(\theta,I)).
\end{equation}
\noindent Recall that $L^{2}({\Bbb T}^{n};{\cal L})$ has a natural
Hilbert basis given by the shifted exponentials, $e_{m}(x)= \exp (
i\langle m + \nu/4, x \rangle ),$ with $m \in {\Bbb Z}^{n}$.
\noindent The next step involves reduction of the order of the
error term on the RHS of (\ref{sec44}). We do this by making
Fourier series expansions of all symbols in terms of the shifted
exponentials, $e_{m}(x)$, and conjugating the expression on the
RHS of (\ref{sec44}) by $U_{1}:= \exp ( i \hbar \,Op_{h}(v_{1}))$
where $v_{1}(\theta,I) \in C^{\infty}({\Bbb T}^{n} \times D_{1})$:
\begin{eqnarray} \label{sec45}
\exp [-i \hbar \,Op_{\hbar}(v_{1})]\circ [ Op_{\hbar}(I_{j}) +
\hbar \, Op_{\hbar}(r_{1})] \circ \exp [ i \hbar
\,Op_{\hbar}(v_{1})] \\ \nonumber
 =_{\Omega}[Id - i \hbar \, Op_{\hbar}(v_{1})]\circ [ Op_{\hbar}(I_{j}) + \hbar
Op_{\hbar}(r_{1})] \circ [ Id + i \hbar \, Op_{\hbar}(v_{1})] +
{\cal O}(\hbar^{2}).
\end{eqnarray}
\noindent Note that, by the pseudodifferential symbolic calculus,
\begin{equation} \label{sec46}
[Id - i \hbar \, Op_{\hbar}(v_{1})]\circ [ Op_{\hbar}(I_{j}) +
\hbar \, Op_{\hbar}(r_{1})] \circ [ Id + i \hbar \,
Op_{\hbar}(v_{1})] + {\cal O}(\hbar^{2}) =_{{\Bbb T}^{n} \times
D_{1}} Op_{\hbar}(I_{j}) + \hbar \,Op_{\hbar}(e_{1j})    + {\cal
O}(\hbar^{2}).
\end{equation}
\noindent Clearly, to get rid of the $Op_{\hbar}(e_{1j})$ term on
the RHS of (\ref{sec46}), we would like to solve the following
transport equations for $v_{1}$:
\begin{equation}  \label{sec47}
e_{1j}(\theta,I) = \{ v_{1}, I_{j} \}(\theta,I) + r_{1}(\theta,I).
\end{equation}
\noindent Keeping in mind that all operators are acting on
$L^{2}({\Bbb T}^{n};{\cal L})$, it follows that these equations
can be rewritten in the form:
\begin{equation}  \label{sec48}
 \frac{\partial}{\partial \theta_{j}} v_{1}(\theta,I) - \frac{\nu}{4}
v_{1}(\theta,I) +r_{1}(\theta,I) = 0\,\,\,\,\,;j=1,...,n.
\end{equation}
\noindent However, in order to solve (\ref{sec48}), we have to
subtract off the zeroth shifted Fourier coefficient,
$\tilde{f}_{j1}(I)$, of $r_{1}(\theta,I)$, and then  match Fourier
coefficients on both sides of (\ref{sec48}). Consequently, if
$$ f_{j1}(q_{1}-b_{1},...,q_{n}-b_{n}) \circ \kappa =
\tilde{f}_{j1}(I_{1},....,I_{n}),$$ \noindent it follows from
(\ref{sec44})-(\ref{sec48}) above that for $j=1,...,n,$
\begin{equation}
U_{1}^{-1}U_{0}^{-1} \circ [ f_{0j}(Q_{1}-b_{1},...,Q_{n}-b_{n}) +
\hbar \, f_{1j}(Q_{1}-b_{1},...,Q_{n}-b_{n}) ] \circ U_{0}U_{1}
=_{{\Bbb T} \times D_{1}} Op_{\hbar}(I_{j}) + {\cal O}(\hbar^{2}).
\end{equation}
\noindent To finish the proof of the lemma, simply apply induction
where the inductive step is    exactly the same as the argument
for the first step given above. \qed

\subsection{Quasimodes}

In this section, we introduce the  quasimodes $
\psi_{b,\hbar}^{(l)}(x) := U^{(l)}_{b}u_{\hbar}(x)$.

Suppose  that $\phi_{b,\hbar}$ is a joint eigenfunction of
$Q_{1},...,Q_{n}$ with $Q_{j} \phi_{b,\hbar} = (b_{j} +
\mu_{j}(\hbar)) \phi_{b,\hbar}$, where $\mu_{j}(\hbar) = o(1).$
Then, as a consequence of the normal form in Lemma (\ref{QBNF}),
$u_{\hbar} =_{\Omega^{(l)}} U^{(l)*} \phi_{b,\hbar}$ is an
admissible, microlocal solution  of the system of equations:
\begin{equation} \label{sec421}
\frac{\hbar}{i} \frac{\partial}{\partial \theta_{j}} u_{\hbar}
=_{\Omega_{0}^{(l)}} m_{j}(\hbar) u_{\hbar}, \,\,\,j=1,2,...,n.
\end{equation}
\noindent where, $m_{j}(\hbar) \in \hbar( {\Bbb Z}+
\frac{\pi}{2}[\nu^{(l)}] )+ {\cal O}(\hbar^{\infty})$ and
$m_{j}(\hbar) = o(1)$ as $\hbar \rightarrow 0$.

\begin{prop} \label{QM1}
The space of admissible microlocal solutions of the system of
equations in (\ref{sec421}) is one-dimensional. Consequently, the
functions
$$u_{\hbar}(\theta)  = \exp [i  m_{1}(\hbar) \theta_{1} + ...+ i m_{n}(\hbar)
\theta_{n})]$$ \noindent span over ${\Bbb C}(\hbar)$ the space of
admissible, microlocal solutions to (\ref{sec421}).
\end{prop}
\noindent{\em Proof:} Let $u_{\hbar}$ be an admissible, microlocal
solution to (\ref{sec421}), by using Fourier transforms and
cutting-off in the frequency variables \cite{CP}, one readily
shows that there exists a global distributional solution,
$v_{\hbar} \in {\cal D }'({\Bbb T}^{n};{\cal L})$, to
(\ref{sec421})  with $v_{\hbar} =_{\Omega_{0}} u_{\hbar}$. But
then, by a simple quadrature, we have that $v_{\hbar}  = c(\hbar)
\,\exp [i ( m_{1}(\hbar)) \theta_{1} + ...+ + m_{n}(\hbar)
\theta_{n})] + {\cal O}(\hbar^{\infty})$ for some $c(\hbar) \in
{\Bbb C}(\hbar)$. \qed

Proposition (\ref{QM1}) has a natural corollary that characterizes
the space of the so-called {\em semiclassical quasimodes} [CV2]
attached to a component torus, $\Omega^{(l)}(b)$:
\begin{cor} \label{QM2}
The space of admissible, microlocal solutions to the system of
equations
$$ Q_{j} \phi_{b,\hbar} =_{\Omega^{(l)}} (b_{j} + \mu_{jk}(\hbar))
\phi_{b,\hbar} + {\cal O}(\hbar^{\infty});\,\,\,j=1,...,n$$ is
spanned over ${\Bbb C}(\hbar)$ by the functions
$$ \psi_{b,\hbar}^{(l)}(x) := U^{(l)}_{b}u_{\hbar}(x),$$
\noindent where, $U^{(l)}_{b}$, is the $\hbar$-Fourier integral
operator constructed in Lemma (\ref{QBNF}) which depends locally
regularly on $b\in B_{reg}$ and $u_{\hbar}(\theta) = \exp [i (
m_{1}(\hbar) \theta_{1} + ...+ m_{n}(\hbar) \theta_{n})]$, with
$m_{j}(\hbar) \in {\Bbb Z}+ \frac{\pi}{2}[\nu^{(l)}] + {\cal
O}(\hbar^{\infty})$.
\end{cor}
\noindent{\em Proof:} The corollary follows immediately from
Proposition (\ref{QM1}) together with the Birkhoff normal form
construction in Lemma (\ref{QBNF}). \qed

\noindent{\bf Remark:} By applying the lemma of stationary phase
in the explicit expressions for the functions
$U^{(l)}_{b}u_{\hbar}(x)$ on a sufficiently small open set $V \in
{\Bbb R}^{n}$, it follows that there exist locally-defined phase
functions, $\phi_{l}(x,\xi) \in C^{\infty}(V \times {\Bbb
R}^{N})$, parametrizing the Lagrangian torus, $\Lambda^{(l)}(b)$,
and classical symbols $b_{l}(x,\xi;\hbar) \sim \sum_{j}
b_{jl}(x,\xi) \hbar^{j}$, such that:
$$ \psi_{b,\hbar}^{(l)}(x) = (2\pi  \hbar)^{-N/2} \int_{{\Bbb R}^{N}}
e^{i\phi_{l}(x,\xi)/\hbar} \, b_{l}(x,\xi;\hbar) \,d\xi.$$
\noindent These oscillatory integrals are precisely the usual
semiclassical Lagrangian distributions (or ``quasimodes'')
attached to Lagrangian tori.

We still have to determine the asymptotic formulas for the
eigenvalues, $ b_{j} + \mu_{j}(\hbar)$ . These formulas are
usually referred to as the (regular) Bohr-Sommerfeld quantization
rules associated with a component Lagrangian torus:

\section{Resonant tori and localization of modes }

 The main results (Lemmas (\ref{RES1})
and (\ref{RES2} ) of this section concern the {\it localization of modes on
tori}.  When the fiber
\begin{equation}\label{FIBRE}  {\cal P}^{-1}(b) = \Lambda^{(1)}(b) \cup \dots \cup  \Lambda^{(m_{cl})}(b), \end{equation}
 where the $\Lambda^{(l)}(b); l=1,...,m$ are
 has multiple components,
it is only true in general that the joint eigenfunctions $\phi_{\lambda}$ with approximately the
joint eigenvalue $\lambda$ concentrate on the union of the tori. The tunnelling
problem which concerns us is the following:
\medskip

\noindent{\bf Problem} How is the mass of $\phi_{\lambda}$ distributed among
the components $\Lambda^{(1)}(b)$ of $ {\cal P}^{-1}(b)$?
\medskip

In the case of a symmetric double well potential $V(x) = V(-x)$ on $\R$, it
is well-known that eigenfunctions are either even or odd, and hence have their
mass equally distributed in both wells. A similar phenomenon occurs in any
integrable system with symmetries. In this case, the wells are resonant, i.e.
the germs of the  Hamiltonian flow at the two tori (or wells) are symplectically
equivalent.

Resonance causes  trouble when comparing sup norms of modes and quasi-modes,
and hence when relating sup norms of modes to  the geometry
of invariant tori. However, this is an exceptional situation: we  prove
that in generic (non-resonant) situations, modes are
well-approximated by individual quasimodes and localize on
individual tori. To deal with the remaining resonant cases, we prove a result (Lemma
(\ref{SUPQM})) which allows us to transfer attention from modes to
quasimodes in the proof of Theorem (\ref{RM}). Namely, we prove
that uniform boundedness of sup norms of modes implies the same
for quasimodes ( Lemma (\ref{SUPQM})).

Much of this  material is implicitly known (it appears in a
slightly different form in \cite{CV2, CP} in the one-dimensional case), but we have not been
able to find the results we need explicitly stated and proved.

\subsubsection{  Bohr-Sommerfeld, Birkhoff normal form and resonant tori }
The notion of resonant tori is based on  the Bohr-Sommerfeld
quantization rules \cite{CV2} for the joint spectrum of
$Q_{1},...,Q_{n}$. In terms of microlocal analysis, these rules
are encoded in the regular quantum Birkhoff `inverse' normal form
construction in Lemma (\ref{QBNF}). Fix a neighbourhood,
$\Omega^{(l)}(b)$, of the torus $\Lambda^{(l)}(b)$ and let
$\phi_{b,\hbar}$ be a microlocal solution of the system of
equations:
\begin{equation} \label{sec431}
Q_{j} \phi_{b,\hbar} =_{\Omega^{(l)}} (b_{j} + \mu_{j}(\hbar))
\phi_{b,\hbar} + {\cal O}(\hbar^{\infty}), \,\,\,\,\mbox{with}
\,\,\,\mu_{j}(\hbar) = o(1)\,\,\,\,;j=1,...,n.
\end{equation}
\noindent Here, by a slight abuse of notation, we have written
$\mu_{j}$ in place of $\mu_{jk}$. In addition, suppose that
\begin{equation} \label{sec432}
\| \phi_{b,\hbar} \|_{\Omega^{(l)}} \geq C
\hbar^{m},\,\,\,\mbox{for some}\,\,m \in {\Bbb Z}.
\end{equation}
\noindent Then, from the `inverse' Birkhoff normal form
construction in Lemma 4.1, it follows that, for $j=1,...,n$,
\begin{equation} \label{sec433}
\sum_{k=0}^{\infty}
f_{jk}^{(l)}(\mu_{1}(\hbar),...,\mu_{n}(\hbar)) \hbar^{k-1} = 2\pi
{\Bbb Z} + \frac{\pi}{2}[\nu^{(l)}] + {\cal O}(\hbar^{\infty}).
\end{equation}
\noindent Let $\gamma_{j}^{(l)}; j=1,...,n$ denote the homology
generators of the torus ${\Bbb T}^{n}$ and $\kappa$ denote the
subprincipal form associated with $Q_{1},...,Q_{n}$ (see, for
instance [GS2]). Then, to first approximation, the equation in
(\ref{sec433}) just the well-known formula:
$$ \frac{1}{\hbar} \int_{\gamma_{j}^{(l)}} \xi dx + \int_{\gamma_{j}^{(l)}}
\kappa +  \frac{\pi}{2}[\nu^{(l)}] = 2\pi {\Bbb Z}  + {\cal
O}(\hbar^{2}).$$ The asymptotic identities in (\ref{sec433}) are
called the {\em Bohr-Sommerfeld (B-S)} quantization rules for the
joint eigenvalues $(b_{1} + \mu_{1}(\hbar),...,b_{n} +
\mu_{n}(\hbar))$ of $(Q_{1},...,Q_{n})$ associated with the
component torus, $\Lambda^{(l)}(b_{1},...,b_{n})$.

Note that in Lemma (\ref{QBNF}), we could have just as easily
derived the Birkhoff normal form as:
\begin{equation}
V^{-1} g(\hbar D_{\theta_{1}},...,\hbar D_{\theta_{n}};\hbar) V
=_{\Omega^{(l)}} Q_{j}-b_{j},
\end{equation}
\noindent where, $g_{j}(x_{1},...,x_{n};\hbar) \sim \sum_{k}
g_{jk}(x_{1},...,x_{n})\hbar^{j}$, the $g_{jk}$'s being the
inverse functions of the $f_{jk}$'s in Lemma (\ref{QBNF}) and
$V:C^{\infty}(\Omega^{(l)}) \rightarrow C^{\infty}({\R}^{n} \times
{\Bbb S}^{1})$ a microlocally unitary $\hbar$-Fourier integral
operator. This way around, we get the (inverse) Bohr-Sommerfeld
quantization rules expressing the
$b_{1}+\mu_{1},...,b_{n}+\mu_{n}$ in terms of the quantum numbers,
$m_{j}(\hbar) \in \hbar({\Bbb Z} + \pi/2 [\nu]) + {\cal
O}(\hbar^{\infty})$:
\begin{equation}
\sum_{k=0}^{\infty} g_{jk}^{(l)}(m_{1}(\hbar),...,m_{n}(\hbar))
\hbar^{k} = \mu_{j}(\hbar) + {\cal O}(\hbar^{\infty}) \,\,\,;
j=1,...,n.
\end{equation}

\noindent In the following, we will need to consider two separate
cases, roughly corresponding to whether or not the functions
$f_{jk}^{(l)}; l=1,...,k$ (or equivalently, the $g_{jk}^{(l)};
l=1,...,k$) are  all the same  for the different component
Lagrangian tori.

We can now define  resonance between tori on the same level of the
moment map.   Let $b\in {\Bbb R}^{n}$ be a fixed regular value of
the moment map ${\cal P}$ and consider the ladder of eigenvalues
$$\Sigma_{b}(\hbar) := \{ (b_{1} + \mu_{1}(\hbar),...,b_{n}+ \mu_{n}(\hbar)) \in
Spec (Q_{1},...,Q_{n}); \, \forall j=1,...,n \, \mu_{j}(\hbar) =
o(1) \}.$$ We say that the component tori,
$\Lambda^{(l)}(b_{1},...,b_{n}); l=1,...,k$, are {\em not in
resonance}  if there exist  eigenvalues
$(\lambda_{1},...,\lambda_{n}) \in \Sigma_{b}(\hbar)$ which
satisfy the Bohr-Sommerfeld (B-S) quantization rules
(\ref{sec433}) associated with the torus,
$\Lambda^{(1)}(b_{1},...,b_{n})$ to ${\cal O}(\hbar^{\infty})$,
but only solve the (B-S) rules associated with the other component
tori, $\Lambda^{(l)}(b_{1},...,b_{n}),$ with $ l \neq 1$ up to
some finite order $m < \infty$. If a pair of tori on the same
level satisfy precisely the same (B-S) quantization conditions,
then they are said to be resonant.

\vspace{2mm}

\subsection{Localization in the non-resonant case}

  The relationship between quasimodes and modes can be fairly complicated when there are  resonances between component Lagrangian tori. However, as we now show, when the generic non-resonance condition is satisfied, the relationship is quite simple.   In the following, we let ${\cal O}_{b}$ denote the
usual order symbol ${\cal O}$ with locally regular dependence on
$b \in B_{reg}$.

\vspace{2mm}

\begin{lem}
 \label{RES1}
Suppose $\phi_{b,\hbar}$ is an $L^{2}$-normalized  joint
eigenfunction of $Q_{1},...,Q_{n}$ with joint eigenvalues
$(\lambda_{1},...,\lambda_{n}) \in \Sigma_{b}(\hbar); j=1,...,n$
 that satisfies the (B-S) rules to infinite order only for the single component
torus, $\Lambda^{(1)}(b_{1},...,b_{n})$ (ie. the component tori,
$\Lambda^{(l)}(b); l=1,...,m$ are not in resonance). Then,
$$ ( Op_{\hbar}(\chi^{(1)}) \phi_{b,\hbar}, Op_{\hbar}(\chi^{(1)})
\phi_{b,\hbar} ) = 1 + {\cal O}_{b}(\hbar^{\infty}).$$ \noindent
Here, $Op_{\hbar}(\chi^{(1)})$ denotes the pseudodifferential
cut-off with symbol $\chi^{(1)}(x,\xi)$ supported in
$\Omega^{(1)}(b)$.
\end{lem}
\noindent {\em Proof:} As we have already shown in Proposition
(\ref{QM1}), the space of microlocal solutions to (\ref{sec421})
is at most one-dimensional. Suppose then that $\phi_{b,\hbar}$ is
a joint eigenfunction of $Q_{1},...,Q_{n}$ with eigenvalues
$(b_{1} + \mu_{1},...,b_{n}+\mu_{n}) \in \Sigma(\hbar)$. In
particular, $\phi_{b,\hbar}$ is then a microlocal solution near
each component torus; that is,
$$ Q_{j} \phi_{b,\hbar} =_{\Omega^{(l)}} (b_{j} + \mu_{j}(\hbar)) \phi_{b,\hbar}
+ {\cal O}_{b}(\hbar^{\infty}).$$ \noindent Consequently, given
Proposition (\ref{QM1}), there are two possibilities: In the first
case, given $\chi^{(l)}(x,\xi) \in C^{\infty}_{0}(\Omega^{(l)})$,
a cutoff supported near the torus, $\Lambda^{(l)}$, we have that
$(Op_{\hbar}(\chi^{(l)}) \phi_{b,\hbar},
Op_{\hbar}(\chi^{(l)})\phi_{b,\hbar}) \geq C \hbar^{m}$ for some
$m \in {\Bbb Z}$. In such a situation, we already know that the
eigenvalues $(b_{1} + \mu_{1},...,b_{n}+ \mu_{n}) \in
\Sigma_{b}(\hbar)$ must satisfy the (B-S) quantization rules
corresponding to the component torus, $\Lambda^{(l)}$, to ${\cal
O}_{b}(\hbar^{\infty})$. The only other possibility is that
$$ ( Op_{\hbar}(\chi^{(l)}) \phi_{b,\hbar}, Op_{\hbar}(\chi^{(l)})
\phi_{b,\hbar}) = {\cal O}_{b}(\hbar^{\infty}).$$ \noindent Given
the definition of non-resonance above, it follows that the latter
scenario holds for $l \neq 1$ and the lemma follows.\qed

As a direct consequence of Lemma \ref{RES1}, it follows that in
the in non-resonant case, $L^{2}$-normalized quasimodes attached
to individual Lagrangian tori agree with actual
joint-eigenfunctions, modulo ${\cal O}(\hbar^{\infty})$ errors in
$L^{2}$. Moreover, we  can now also prove the following
correspondence principle for the case of ${\Bbb R}^{n}$ actions:
\begin{prop}  \label{RES2}
Let $q(x,\xi) \in C^{\infty}_{0}(T^{*}M)$ and suppose that the
component tori, $\Lambda^{(l)}(b); l=1,...,m$ are not in
resonance. Then, for all  $L^{2}$-normalized joint eigenfunctions,
$\phi_{b,\hbar}$, with joint eigenvalues $(b_{1}+\mu_{1}(\hbar),
... , b_{n}+ \mu_{n}(\hbar)) \in \Sigma_{b}(\hbar)$ with $
\mu_{j}(\hbar) = {\cal O}(\hbar)$, there exists a component torus,
$\Lambda^{(1)}(b)$, such that:
$$ ( Op_{\hbar}(q) \phi_{b,\hbar}, \phi_{b,\hbar}) = (2\pi)^{-n}
\int_{\Lambda^{(1)}(b)} q \,d\mu_{b} + {\cal O}_{b}(\hbar).$$
\noindent Here, $d\mu_{b}$ denotes Liouville measure on
$\Lambda^{(1)}(b)$.
\end{prop}
\noindent {\bf Proof}: The proof is  similar to the case of torus
actions \cite{Z1}. For the sake of completeness, we sketch the
argument. First, in light of Lemma (\ref{RES1}),
\begin{equation} \label{sec451}
 ( Op_{\hbar}(q) \phi_{b,\hbar}, \phi_{b,\hbar}) = ( Op_{\hbar}(q)\, \circ
Op_{\hbar}(\chi^{(1)}) \phi_{b,\hbar}, Op_{\hbar}(\chi^{(1)})
\phi_{b,\hbar}) + {\cal O}(\hbar^{\infty}).
\end{equation}
\noindent A standard application of averaging, combined with  the
semiclassical Egorov theorem and the normal form in Lemma  gives:
\begin{equation} \label{sec452}
 ( Op_{\hbar}(q)\, \chi^{(1)} \phi_{b,\hbar}, \chi^{(1)} \phi_{b,\hbar}) =
(2\pi)^{-n} \int_{\Lambda^{(1)}} q \,d\mu  + e(\hbar) + {\cal
O}(\hbar),
\end{equation}
\noindent where,
 $$ e(\hbar) = ( Op_{\hbar}(r) u_{\hbar}, u_{\hbar} )$$
\noindent and $r \in C^{\infty}_{0} ({\Bbb T}^{n} \times D_{1})$
with $r(\theta,I) = {\cal O}(|I|)$ and $u_{\hbar}(\theta) = \exp
[i(m_{1} \theta_{1} + ...+ m_{n}\theta_{n})]$. An integration by
parts in  the $I_{1},...,I_{n}$ variables  shows that:
$$( Op_{\hbar}(r) u_{\hbar}, u_{\hbar} ) = {\cal O}(\hbar).$$
\noindent The proposition follows.\qed

\subsection{Localization in the resonant case}

In resonant cases,  modes do not
necessarily localize on individual tori.  We introduce the
language of `quantum limit measures' to describe localization in
more general cases.

\subsubsection{Quantum limits of modes and quasimodes}

Let us denote by $M_I$ the set of measures on $S^*M$ which are
invariant under the Hamiltonian $\R^n$-action. Among such measures
are the orbital averaging measures $\mu_z$, defined by
$\int_{S^*M} f d \mu_z = \lim_{T \to \infty} \frac{1}{T^n}
\int_{\max |t_j| \leq T}  f(\Phi^t(z)) dt.$ In the case of compact
(torus) orbits, one gets the normalized (probability) Lebesgue
measure on the torus.

By the `quantum limit' measures ${\mathcal Q} $ of the quantum
integrable system, we mean the set of weak* limits of the measures
$d \Phi_{\lambda}$ defined by
\begin{equation} \langle Op(a) \phi_{\lambda}, \phi_{\lambda}\rangle =
\int_{S^*M} a d \Phi_{\lambda}, \end{equation} where $Op(a)$ is
the pseudodifferential operator (in some fixed quantization) with
complete symbol $a$. It is an easy consequence of Egorov's theorem
that ${\mathcal Q} \subset M_I.$ For background,
 terminology and references in a  context closely related to this one,  we refer
to \cite{TZ}.

When $b$ is a regular value of the moment map ${\cal P}$,
$ {\cal P}^{-1}(b)$ is given by (\ref{FIBRE}), where the $\Lambda^{(l)}(b); l=1,...,m$ are
$n$-dimensional Lagrangian tori. Let $d\mu_{\Lambda^{(j)}(b)}$
denote  the normalize Lesbegue measures on such tori.  If
$\{\phi_{\lambda_j}\}$ is a sequence of eigenfunctions with
$\frac{\lambda_j}{|\lambda_j|}\} \to b,$ then we have (see
\cite{T2})
$$d \Phi_{\lambda_j} \to \sum_{j = 1}^m c_j
d\mu_{\Lambda^{(j)}(b)},\;\;\;\mbox{for some}\;\;\; c_j \geq 0,
\sum c_j = 1.$$ It is difficult to determine the coefficients in
general. As we have seen,
 in the
generic (`non-resonant') case, one can find for each component
torus a sequence of eigenfunctions for which $c_j = 0, j \geq 2.$

Quite analogous to the quantum limit measures are the analogous
measures defined for quasi-modes associated to regular tori. We
thus define the set ${\mathcal Q}_{qm} \subset M_I$ of `quasimode
limit measures' to be the
 weak* limits of the analogous measures $d \Psi_{\lambda}$ corresponding to
quasimodes $\{\psi_{\lambda}\}$ associated to regular tori. The
following lemmas will be crucial for the proof of the Theorem in
the ${\Bbb R}^{n}$ case:

\begin{lem} \label{QML} For
any component Lagrangian torus $\Lambda^{(l)}(b)$,  there exists a
sequence $\psi_{b,\hbar}^{(l)}$  of quasimodes such that
\begin{equation}
(Op_{\hbar}(q) \psi^{(l)}_{b,\hbar}, \psi^{(l)}_{b,\hbar}) =
(2\pi)^{-n} \int_{ \Lambda^{(l)}} q \,\,d\mu_{b} + {\cal
O}_{b}(\hbar),
\end{equation}
i.e. that ${\mathcal Q}_{qm}$ contains all normalized Lebesgue
measures along compact Lagrangean orbits. \end{lem}

\noindent {\bf Proof:} Since the quasimodes $\psi_{b,\hbar}^{(l)}$
are, by construction, microlocally concentrated on a particular
component torus, $\Lambda^{(l)}(b),$ Lemma (\ref{QML}) is proved
in the same way as Proposition (\ref{RES2}) above. \qed

\section{Proof of Theorem (\ref{RM})}

\subsection{Pointwise bounds on quasimodes}

 In order to prove Theorems \ref{RM} in the possibly resonant case, we will need to establish some less direct relations between modes and quasimodes. The next Lemma shows that uniform bounded of sup-norms of   modes implies the same
for quasimodes:

\begin{lem} \label{SUPQM} Assume  $\| \phi_{\lambda} \|_{\infty} \leq A$ for
all $\lambda \in Spec(P_{1},...,P_{n})$. Then for all regular $b
\in {\Bbb R}^{n}$, the $L^2$-normalized quasimodes satisfy:
$$ \limsup_{h \to 0} \| \psi^{(l)}_{b, h} \|_{\infty} \leq A \;\;m_{cl}(b).$$
\end{lem}

\bigskip

\subsubsection{Proof of Lemma (\ref{SUPQM})}
\noindent {\bf Proof:}

 Given the component tori $\Lambda^{(l)}(b)$ where $l=1,...,m$, recall that we have  denoted by
$\psi_{b,\hbar}^{(l)}$ the quasimodes in Corollary (\ref{QM2})
with approximate eigenvalues $(b_{1} + \mu_{1},...,b_{n}+ \mu_{n})
\in \Sigma_{b}(\hbar)$ which are microlocally concentrated near
the component torus, $\Lambda^{(l)}(b)$, . By making a
Gram-Schmidt orthogonalization, we can assume that these functions
are orthonormal in $L^{2}$. Recall, they satisfy for each
$j=1,...,n$:
$$ Q_{j} \psi_{b,\hbar}^{(l)} = (b_{j} + \mu_{j}(\hbar)) \psi_{b,\hbar}^{(l)} +
{\cal O}_{b}(\hbar^{\infty})$$ \noindent in $L^{2}$. Given the
(B-S) rules, we know that for any component torus
$\Lambda^{(l)}(b)$, there exist $\epsilon_{b}^{(l)} >0$ such that
for $\hbar$ sufficiently small, the gaps in the (B-S) lattice
spectrum are greater than $\epsilon_{b}^{(l)} \, \hbar$. Define
$$ \epsilon_{b} := \min_{l \leq m_{cl}(b)} \epsilon_{b}^{(l)},$$
\noindent and consider the cube,
$$ C_{b,\hbar}:= \prod_{j=1}^{n} \,[ -\epsilon_{b} \,\hbar + b_{j} +
\mu_{j}(\hbar), \, \epsilon_{b} \,\hbar + b_{j} + \mu_{j}(\hbar)
\,].$$

 Let ${\bf Q} := (Q_{1},...,Q_{n})$ act on $L^{2}(M) \oplus \cdot
\cdot \cdot \oplus L^{2}(M)$ componentwise. Then, we denote by
${\cal Q}^{(l)}_{b,\hbar}$
 the (at most) one-dimensional vector space generated by the quasimodes attached to $\Lambda^{(l)}(b)$ with
quasieigenvalues in $C_{b,\hbar}$. In addition, we define
\begin{equation} \label{sec441}
{\cal Q}_{b,\hbar} := \bigoplus_{l=1}^{m_{cl}(b)} {\cal
Q}^{(l)}_{b,\hbar}.
\end{equation}
\noindent So, put another way, ${\cal Q}_{b,\hbar}$ is the vector
space of functions generated by quasimodes $\psi_{b,\hbar}^{(l)};
l=1,...,m_{cl}$ from all of the component tori with
quasieigenvalues in the cube $C_{b,\hbar}$. Similarily, we denote
by   ${\cal R}_{b,\hbar}$ the vector space generated by joint
eigenfunctions of ${\bf Q}$ associated with joint eigenvalues of
${\bf Q}$ lying in $ C_{b,\hbar}$.  As a consequence of the (B-S)
rules, we know that there are gaps in the joint spectrum of the
${\bf Q}$: that is, there exists $C_{b}>0$ such that for all
$\hbar \in (0,\hbar_{0}(b)]$,
$$dist \, ( \, Spec({\bf Q})- Spec({\bf Q} \cap C_{b,\hbar} )\, , \, Spec({\bf
Q} \, ) \, ) \, \geq \frac{1}{C_{b}}\hbar.$$
 Let $E_{{\cal R}_{b,\hbar}}$ denote the
spectral projector onto  the eigenvalues of ${\bf Q}$ lying in
${\cal R}_{b,\hbar}$ and consider a fixed quasimode,
$\psi_{b,\hbar} \in {\cal Q}_{b,\hbar}$. It follows that:
\begin{equation} \label{sec442}
{\bf Q}  (\psi_{b,\hbar} - E_{{\cal R}_{b,\hbar}}\psi_{b,\hbar})
= {\cal O}_{b}(\hbar^{\infty}).
\end{equation}
\noindent So, by the spectral theorem,
\begin{equation} \label{sec443}
\| \psi_{b,\hbar} - \sum_{\phi_{b,\hbar} \in R_{b,\hbar}}
(\psi_{b,\hbar}, \phi_{b,\hbar}) \phi_{b,\hbar} \|_{(0)} = {\cal
O}_{b}(\hbar^{\infty}).
\end{equation}
\noindent Conversely, given $\phi_{b,\hbar} \in {\cal
R}_{b,\hbar}$, it follows from the microlocal characterization
result in Proposition (\ref{QM1}) that:
\begin{equation} \label{sec444}
\| \phi_{b,\hbar} - \sum_{\psi_{b,\hbar} \in {\cal Q}_{b,\hbar}}
(\phi_{b,\hbar}, \psi_{b,\hbar}) \psi_{b,\hbar} \|_{(0)} = {\cal
O}_{b}(\hbar^{\infty}).
\end{equation}
\noindent Let $d({\cal Q}_{b,\hbar},{\cal R}_{b,\hbar}) = \|
\Pi_{{\cal Q}_{b,\hbar}} - \Pi_{{\cal Q}_{b,\hbar}} \Pi_{{\cal
R}_{b,\hbar}} \|_{(0)}$ denote the non-symmetric distance between
the vector spaces ${\cal Q}_{b,\hbar}$ and ${\cal R}_{b,\hbar}$.
It follows from (\ref{sec443}) and (\ref{sec444}) that, in
particular, $d({\cal Q}_{b,\hbar},{\cal R}_{b,\hbar}) <1$ and
$d({\cal R}_{b,\hbar},{\cal Q}_{b,\hbar}) <1$. Consequently [H],
${\cal Q}_{b,\hbar} \cong {\cal R}_{b,\hbar}$. Since ${\cal
Q}_{b,\hbar}$ has dimension uniformly bounded independent of
$\hbar$, so does ${\cal R}_{b,\hbar}$. In particular, we have that
\begin{equation} \label{dim}
\dim \, {\cal R}_{b,\hbar} \, \leq \, m_{cl}(b).
\end{equation}

Given the estimates in (\ref{sec442})-(\ref{dim}), we have by the
Garding and Sobolev inequalities,
\begin{equation}  \| \psi^{(l)}_{b, h} \|_{\infty} \leq \| E_{R_{b,h}}
\psi^{(1)}_{b, h}\|_{\infty} + C \| (I -
\Delta)^{s}(\psi^{(l)}_{b,\hbar} - E_{R_{b,h}} \psi^{(1)}_{b,
h})\|_{(0)}.
\end{equation}
The lemma follows from the following two estimates:
\begin{equation} \begin{array}{l}  \| E_{R_{b,h}} \psi^{(1)}_{b, h}\|_{\infty} =
\|\sum_{j = 1}^{m(b)} \langle \psi^{(1)}_{b, h}, \phi_{b,h}^{(j)}
\rangle \phi_{b,h}^{(j)} \|_{\infty}
\leq A \;\; m_{cl}(b)\\ \\
\| (I - \Delta)^{s}(\psi^{(l)}_{b,\hbar}  - E_{R_{b,h}}
\psi^{(1)}_{b, h})\|_{(0)} = {\cal O}_b(\hbar^{\infty})
\end{array} \end{equation}

\noindent \qed

\subsection{Completion of the proof}

The completion of the proof is now exactly as in \cite{TZ}, so we only
briefly recall it. We first observe:

\begin{lem}\label{UB} Let $\{ \pi_* d \mu_{\Lambda}\}$ denote the set of
projections to $M$ of normalized Lebesgue measures on compact
Lagrangean tori $\Lambda \subset S^*M.$ Then the family is
uniformly bounded as linear functionals on $L^{1}(M).$ \end{lem}

\begin{proof} By Lemma (\ref{QML}), each $ d
\mu_{\Lambda}$ is a weak limit of quasi-mode measures
corresponding to a sequence, say  $\{\psi_{\Lambda, h}\}$,  of
quasimodes. By Lemma (\ref{SUPQM}) we have, for $V \in C(M)$,
\begin{equation} \begin{array}{l} |\int_{M} V \pi_* d\mu_{\Lambda}| = \lim_{h
\to 0} |\langle V \psi_{\Lambda, h}, \psi_{\Lambda, h} \rangle|
\leq ||V||_{L^1} \limsup_{h \to 0} ||\psi_{\Lambda, h}||^{2}_{\infty}\\ \\
 \leq A^{2} \;\; m_{cl}({\mathcal P}(\Lambda))^{2}\;\;  ||V||_{L^1}. \end{array}
\end{equation} \end{proof}

\begin{cor} The invariant Lagrangean tori project without singularities to
the base $M$. \end{cor}

\begin{proof}  Indeed, we have $ \pi_* d \mu_{\Lambda} = f_{\lambda} dvol$ for
some $f_{\lambda} \in L^{\infty}(M).$ But it is straightforward to see that
if $\pi|_{\Lambda}$ has singularities, then $f_{\lambda}$ will blow up along them.
\end{proof}

The geometric problem then arises, which  completely integrable systems
have  the property that all invariant tori project without
singularities to the base $M$? It follows immediately that $M$ must be a torus,
and that no leaf of the Liouville foliation is singular. By a theorem of Mane,
it then follows that the metric has no conjugate points. But the only metrics
on a torus without conjugate points are the flat ones (Burago-Ivanov). For
more details, we refer to \cite{TZ}.

\section{\label{EXAMPLES} Examples of quantum completely integrable systems}

To illustrate some of the variety of quantum completely integrable
systems, we  describe the Liouville foliations of two infinite
dimensional families of two-dimensional examples: surfaces of
revolution and Liouville tori.

\subsection{Spheres of revolution}

Suppose that $g$ is a metric on $S^2$ for which there exists
 an effective action of $S^1$ by isometries of
$(S^2, g)$.  The two fixed points will be denoted $N,S$ and $(r,
\theta)$ will denote geodesic polar coordinates centered at $N$,
with $\theta = 0$ some fixed meridian $\gamma_M$ from $N$ to $S$.
The metric may then be written in the form $g = dr^2 + a(r)^2
d\theta^2$ where $a: [0,L] \rightarrow \R^+$ is defined by $a(r) =
\frac{1}{2\pi} |S_r(N)|$, with $|S_r(N)|$  the length of
 the distance circle of radius $r$ centered at $N$.  For any smooth surface
of revolution, the function $a$ satisfies $a^{(2p)}(0) =
a^{(2p)}(L) = 0, a'(0) = 1, a'(L) = -1$ and two such surfaces
$M_1, M_2$ are isometric if and only if $L_1 = L_2$ and $a_1(r) =
a_2(r)$ or $a_1(r) = a_2(L-r).$

Complete integrability of $H = |\xi|_g$ (i.e. of the geodesic
flow) is classical, and follows from the existence of the Clairaut
integral $p_{\theta}(v):= \langle v, \frac{\partial}{\partial
\theta}\rangle$.  Since
 the Poisson bracket $\{p_{\theta}, |\xi|_g\} = 0$, the geodesics
are constrained to lie on the level sets of $p_{\theta}$; and
since both $|\xi|_g$ and $p_{\theta}$ are homogeneous of degree
one, the behaviour of the geodesic flow is determined by its
restriction to $S^*_gS^2 = \{|\xi|_g = 1\}$.
 The moment map is given by
$$P = (|\xi|_g, p_{\theta}): T^*S^2 \rightarrow \R^2.$$

To describe the inverse images, i.e. the invariant tori,  we need
to specify the function $a(r).$  We will assume it  is a Morse
function on $[0,L]$ with critical points at $\{m_{k} \}\subset
(0,L)$.

\subsubsection{Case 1: generic case}

 For simplicity and with no real loss of
generality, let us assume that there are two local maxima $M_1,
M_2$ with $a(M_2) = C_2 > a(M_1)= C_1$ and one
 local minimum $m_2$ with $a(m_2) = c_2.$
We note that the canonical involution $(x, \xi) \rightarrow (x,
-\xi)$ preserves ${\cal P}^{-1}(p_1, p_2)$, and that only the
meridian torus ${\cal P}^{-1}(p_1, 0)$
 is invariant; otherwise there are always at least
two components on each level.  By homogeneity it suffices to
consider the level sets ${\cal P}^{-1}(1, c)$ in the unit tangent
bundle. We will denote the  components  by $M^{(\pm, k)}_{c}$
where $\pm$ denotes some ordering of each pair of tori
interchanged by the involution. The following table summarizes the
invariant tori, their projections to $M$ and the type of
singularity of the projection.
\bigskip

\begin{tabular}{llll} Range of $c$ &  components $M^{(\pm, k)}_{c}$   &
projections to $M$ & singularities \\
(0) $c > C_2$ or $c < 0$ & none & -- & --  \\
(a) $c = C_2$ & 2 circles & Project to closed geodesics & singular leaf \\
(b) $C_1 < c < C_2$ & 2  torii  & 2 horizontal  annuli  & fold \\
(c) $ c_1 < c < C_1$ & 4  tori & each projects diffeo to $M$ & none\\
(d) $c = c_{1 }$ & 4 cylinders  & each projects over all of $M$ & none \\
(e) $0 < c < c_1$ & 2 tori & projects to 2 horizontal annuli & fold\\
(f) $c = 0$ & 1 (meridian) torus & covers all of $M$& blow-down.\\
\end{tabular}
\bigskip

\subsubsection{ Simple surfaces of revolution}

Simple surfaces of revolution are examples of toric integrable
systems \cite{CV1}  They are the special surfaces where:
\medskip

\noindent(i) $a$ has precisely one non-degenerate critical point
$r_o \in (0, L)$, with $a''(r_o) < 0$,
corresponding to an `equatorial geodesic' $\gamma_E$;\\
\noindent(i)  the (non-linear) Poincare map ${\cal P}_{\gamma_E}$
for
$\gamma_E$ is of twist type (cf. \S 1). \\

It is proved in \cite{CV1} that the geodesic flow of a simple
surface of revolution is    is  toric integrable.

With the extra assumption on $a$, the level sets are compact and
the only critical level is that of the equatorial geodesics
$\gamma_E^{\pm} \subset S^*_gS^2$ (traversed with either
orientation).  The other level sets  consist of two-dimensional
torii.

As above, let
$$P = (|\xi|_g, p_{\theta}): T^*S^2 \rightarrow B:= \{(b_1, b_2) :
|b_2| \leq a(r_o)b_1\} \subset \R \times\R^+ \leqno(1.1.2)$$  be
the moment map of the Hamiltonian $\R^2$-action defined by the
geodesic flow and by rotation.  The singular set of $P$ is the
closed conic set $Z:= \{(r_o, \theta, 0, p_{\theta}): \theta \in
[0, 2 \pi), p_{\theta} \in \R\}$, i.e. $Z$ is the cone thru the
equatorial geodesic (in either orientation). The image of $Z$ is
the boundary of $B$; the map $P|_{T^*S_gS^2 - Z}$ is a trivial
$S^1 \times S^1$ bundle over the open convex cone $B_o$ (the
interior of $B$).
 For each $b\in B_o$ , let
$H_1(F_b, \Z)$ denote the homology of the fiber $F_b := P^{-1}
(b).$ This lattice bundle is trivial since $B$ is contractible, so
there exists a smoothly varying homology basis $\{\gamma_1(b),
\gamma_2(b)\} \in H_1(F_b, \Z)$ which equals the unit cocircle
$S^*_NS^2$ together with the fixed closed meridian $\gamma_M$ when
$b$ is on the center line $\R^+ \cdot (1,0)$. The action variables
are given by \cite{CV1}, \S 6,
$$I_1 (b) = \int_{\gamma_1(b)} \xi dx = p_{\theta},\;\;\;\;\;\;I_2 (b) =
\int_{\gamma_2(b)} \xi dx =\frac{1}{\pi}
\int_{r_{-}(b)}^{r_{+}(b)} \sqrt{b_1^2 - \frac{b_2^2}{a(r)^2}} dr
+ |b_2|\leqno(1.1.3)$$ where $r_{\pm}(b)$ are the extremal values
of $r$ on the annulus $\pi(F_b)$ (with $\pi : S^*S^2 \rightarrow
S^2$ the standard projection). On the torus of meridians in
$S^*S^2$, the value of $I_2$ equals $\frac{L}{\pi}$ and it equals
one on the equatorial geodesic. So extended, $I_1, I_2$ are smooth
homogeneous functions of degree 1 on $T^*S^2$, and   generate
$2\pi$-periodic Hamilton flows. It follows that the pair ${\cal
I}:=(I_1, I_2)$ generates a global Hamiltonian torus ($S^1 \times
S^1$)-action commuting with the geodesic flow.  We refer to it as
the `standard' Hamiltonian $T^2$ action on $T^*S^2-0$. Up to
symplectic equivalence, it is unique.

The singular set of ${\cal I}$ equals ${\cal Z}:= \{I_2 =\pm
p_{\theta}\}$, corresponding to the equatorial geodesics.
 The map
$${\cal I}: T^*S^2 - {\cal Z} \rightarrow \Gamma_o:= \{(x,y) \in \R \times
\R^+ : |x| < y\}\leqno(1.1.4)$$ is a trivial torus fibration.
Henceforth we often write $T_I$ for the torus ${\cal I}^{-1} (I)$
with $I \in \Gamma_o$ and let $\Gamma = {\cal Cl}\Gamma_o$ be the
closure of $\Gamma$ as a convex cone.
 The symplectically dual angle variables
$(\phi_1 = \theta, \phi_2)$ then give, by definition,  the flow
times (mod $2\pi$)
 along the orbits of $(I_1, I_2)$ from a fixed point on $F_b$, which we may
take to be the unique point lying above the intersection of the
equator and the fixed meridian on $F_b$ with the geodesic pointing
into the northern hemisphere.

The following table summarizes the invariant tori and their
projections
\bigskip

\begin{tabular}{llll} Range of $c$ &  components $M^{(\pm, k)}_{c}$   &
projections to $M$ & singularity\\
(0) $c > C$ or $c < 0$ & none & -- & --  \\
(a) $c = C$ & 2 circles & Project to closed geodesics & singular leaf\\
(b) $0 < c < C$ & 2  torii  & 2 horizontal  annuli & fold \\
(c) $c = 0$ & 1 (meridian) torus & covers all of $M$ & blow-down.\\
\end{tabular}
\bigskip

\subsection{ Tori of revolution}  Here, we consider a metric $g$ on $T^2$
which possesses a free effective action $R_{\theta}$ of $S^1$ by
isometries.   We denote the generator of this action by $Z$.
Associated to $Z$ is the Clairaut integral $p_Z : T^*T^2
\rightarrow \R$ defined by $p_Z (\nu) = \langle \nu, Z \rangle.$
The $S^1$ action on $T^2$ lifts to $T^*(T^2)$ as the
 Hamiltonian flow $\phi^t$ of $p_{\theta}.$  Since $\{|\xi|_g, p_Z \} = 0,$ we
have a Hamiltonian $G= \R^1 \times S^1$ action $\Phi^t  = G^t
\times R_{\theta}$ on $T^*(T^2)$ with moment map ${\cal P} =
(|\xi|_g, p_Z).$

Consider now the level set
$$M = {\cal P}^{-1}(1, 0) = \{(x, \xi): |\xi|_g = 1, \langle \xi,
\partial/\partial \theta \rangle = 0\}$$
It is an affine Lagrangean submanifold of $T^*(T^2)$ and since  it
is  compact it must be a Lagrangean torus.  It follows that the
stabilizer $G_m$ of a point of $m \in M$ must be a subgroup of the
form $L \Z \times \{id\}$.  In particular, the $G^t$-orbit of each
$m \in M$ must be closed.

 As above, let us restrict the natural projection $\pi: S^*_g T^2
\rightarrow T^2$ to $M$, to get $\pi_M : M \rightarrow M.$
Clearly, $\pi_M$ is equivariant with respect to the $S^1$ action
$R_{\theta}$ on $T^2$ and its lift.   It sends the $G^t$-orbits of
$m \in M$ to closed geodesics of $T^2$ and since $p_Z = 0$ on $M$,
these geodesics are orthogonal to the $R_{\theta}$-orbits, which
we will call `parallels'.  Further, we note that the tangent space
$T_mM$ is spanned by the
 Hamilton vector fields $\Xi_{|\xi|_g}$, resp. $\Xi_{p_Z}$.  Since $\pi_{M
*}$ sends these vectors to orthogonal vectors, $\pi_M$ is a
submersion and hence a covering map.  Its degree is one because of
the equivariance and the fact that it sends $G^t$-orbits to closed
geodesics (this also follows from Lalonde's theorem cited above).
Hence it is a diffeomorphism.  We will refer to the geodesics in
$M$ as the meridians.

Fixing a point $m_o \in M$ and $x_o = \pi_M(m_o) \in T^2$ we get a
coordinate system on $T^2$ with coordinates $(r, \theta): x = (r,
\theta)  \Leftrightarrow \Phi^{(r, \theta)}(m_o) = \pi_M^{-1}(x).$
Since the meridians are orthogonal to the `parallels',  this forms
an orthogonal coordinate system for $g$.  Hence $g$ has the form:
$g = dr^2  + a(r)^2 d \theta^2.$ \bigskip

\subsubsection{ Liouville torii \cite{B.K.S}\cite{KMS}}  A Riemannian torus $(M, g)$ of
 dimension 2 is called a {\it Liouville torus} if the metric has the form
$$g = [U_1(x_1) - U_2(x_2)](dx_1^2 + dx_2^2) $$
in the coordinates $(x_1, x_2)$ on $[0,1]\times [0,1].$  The
functions $U_j$ are assumed to be periodic of period 1 and to
satisfy $U_1(x_2) - U_2(x_2) > 0$.  The geodesic flow is the
Hamilton flow on $T^*M$ generated by
$$H(x,\xi)^2 = \frac{1}{[U_1(x_1) - U_2(x_2)}(\xi_1^2 + \xi_2^2).$$
It Poisson commutes with the first integral
$$S(x,\xi) = \frac{U_2(x_2)}{[U_1(x_1) - U_2(x_2)}\xi_1^2 +
\frac{U_1(x_1)}{[U_1(x_1) - U_2(x_2)}\xi_2^2.$$ We then have the
homogeneous moment map
$${\cal P} = (p_1, p_2) = (H, \frac{S}{H}): T^*M - 0 \rightarrow \R^2.$$
We note that the canonical involution $(x, \xi) \rightarrow (x,
-\xi)$ preserves ${\cal P}^{-1}(p_1, p_2)$, and that no torus is
invariant; hence there are always at least two components on each
level.  As above, it suffices to consider the level sets in the
unit tangent bundle, which we continue to denote by
 $M^{(\pm, k)}_{c}$.

To describe the inverse images, i.e. the invariant tori,  we need
to specify the functions $U_j(x_j).$  We will assume that $U_j$ is
a Morse function on $[0,1]$ with critical points at $\{c_{jk}
\}\subset (0,1)$ with critical values  $\{v_{jk} \}.$   We briefly
review the description of the tori under the assumption that each
$U_j$ has precisely one  maximum $c_{j+}$ and one minimum $c_{j-}$
and that the critical points do not coincide (see \cite{KMS} for a
detailed description).  We thus have $c_{1+} > c_{11} > c_{2 +} >
c_{2 -}.$

\begin{tabular}{llll} Range of $c$ &  components $M^{(\pm, k)}_{c}$   &
projections to $M$ & singularity \\
(0) $c > c_{1+}$ or $c < c_{2 -}$ & none & -- & -- \\
(a) $c = c_{1+}$ & 2 circles & Project to closed geodesics & singular leaf\\
(b) $c_{1 -} < c < c_{1 +}$ & 2  torii  & 2 horizontal  annuli & folds  \\
(c) $c = c_{1-}$ & 4 cylinders  & each projects over all of $M$  & none\\
 (d) $c_{2 + } < c < c_{1 -}$ & 4  tori & each projects diffeo to $M$ & none \\
(a') $c = c_{2 +}$ & analogous to (a) & 2 circles & singular leaves \\
(b') $c_{2 - } < c < c_{2 +}$ & analogous to (b) & 2 vertical annuli & folds\\
(c') $c = c_{2-}$ & analogous to (c) & &
\end{tabular}
\bigskip

\subsubsection{ Tori of revolution: Quantum theory}  Consider  the metric on
$\R^2/\Z^2$
$$g = a(x)(dx^2 + d\xi^2)$$
where $a$ is a periodic function.  Translations in the $\xi$
variable are then isometries so this is a torus of revolution.
The Laplacian
$$\Delta = \frac{1}{a} \Delta_0$$
commutes with $\partial/\partial \xi$ and hence there is an ONB of
 joint eigenfunctions of the form
$$\phi_{N, \lambda} = e^{2\pi i N \xi} \psi_{N, \lambda}(x).$$
The eigenvalue problem is separable in these variables and gives
$$\psi'' - 4\pi^2 N^2 \psi = - \lambda a \psi.$$
Bourgain [Bour] observes that for $a = 1 - \tau sin^2 \pi x$ and
for sufficiently large $N$ there are eigenfunctions satisfying
$$||\phi_{\lambda}||_{\infty} \geq \lambda^{\frac{1}{8}},
\;\;\;\;\;\;\;||\phi_{\lambda}||_{L^2} = 1.$$

These eigenfunctions concentrate near the projection of a stable
geodesic given by the parallel $x=0$. Note that Theorem 0.2 gives
in this case a lower bound of $\|\phi_{\lambda}\|_{\infty} \geq
C(\epsilon) \lambda^{\frac{1}{8}-\epsilon}$ for any $\epsilon >0$.
Since $da(x) = - 2 \pi \sin \pi x \,\cos \pi x$ and since $x \in
[0,1]$ there are two critical circles: at $x = 0$ and at $x =
\half$, a maximal parallel and a minimal parallel (which are
geodesics).  There are four embedded cylinders consisting of
geodesics which spiral in towards the latter circle.

\subsubsection{ Liouville tori: Quantum theory \cite{B.K.S}\cite{KMS}}
The Laplace Beltrami operator on a Liouville torus, $Q$, with
metric $g = [U_1(x_1) - U_2(x_2)](dx_1^2 + dx_2^2)$ is just:
$$ \Delta = [U_{1}(x_{1})- U_{2}(x_{2})]^{-1} \left( \frac{\partial^{2}}{ \partial x_{1}^{2}} +  \frac{\partial^{2}}{ \partial x_{2}^{2}} \right).$$
The commuting operator with principal symbol $S(x,\xi)$ in 5.4 is
just:

$$S = [U_1(x_1) - U_2(x_2)]^{-1} \left(  U_2(x_2) \frac{\partial^{2}}{ \partial x_{1}^{2}}  +
U_1(x_1) \frac{\partial^{2}}{ \partial x_{2}^{2}} \right).$$

\noindent The joint eigenfunctions of $\Delta$ and $S$ can be
described by two coupled O.D.E.(see K.M.S.). Moreover, in the case
where $c=c_{1+}$ or $c=c_{2+}$ the two circles in the joint level
set are stable and elliptic and project to a single geodesic on
$Q$. In the case $c=c_{1+}$ this projected geodesic is identified
with the line given by $x_{1}=M_{1}$ where $U_{1}(x_{1})$ attains
its non-degenerate maximum at this point, whereas, when $c=c_{2+}$
this geodesic gets identified with the line $x_{2}=M_{2}$ where
$U_{2}(x_{2})$ attains its maximum at $M_{2}$. There are
eigenfunctions attached to these levels which satisfy $\|
\phi_{\lambda} \|_{\infty} \geq \lambda^{\frac{1}{8}}$ and again
this is consistent with the prediction given by Theorem 0.2. Note
also that there exist hyperbolic geodesics on the level sets given
by $c=c_{1-}$ and $c=c_{2-}$ which are limit sets for the
geodesics on these cylinders. These are the lines given by
$x_{1}=m_{1}$ and $x_{2}=m_{2}$ respectively, where $m_{1}$ and
$m_{2}$ are the points where $U_{1}$ and $U_{2}$ achieve their
respective minima. In this case, [T1], one can show that there are
joint eigenfunctions concentrating near the projections of these
hyperbolic orbits and satisfy $\| \phi_{\lambda} \|_{\infty} \geq
\lambda^{\frac{1}{8}} (\log \lambda)^{- \frac{1}{2}}$, slightly
better than the $\| \phi_{\lambda} \|_{\infty} \geq C(\epsilon)
\lambda^{\frac{1}{8}-\epsilon}$ given by Theorem 0.2.

\subsection{Other examples}
 Other classical examples include the Euler top
 and geodesic on asymmetric ellipsoids in any dimension.
 These examples are also quantum completely integrable. Tri-axial
 ellipsoids have special umbilic points where all directions are
 geodesic loop directions. Hence they are candidates for maximal
 eigenfunction growth. The loops form a Lagrangean manifold with
 boundary which has a blow-down projection over the umbilic
 points. We do not believe however that the eigenfunctions
 actually do have maximal sup norm growth.

\section{Further discussion and open problems}

In \cite{TZ}, we list a number of (still) open problems in this
area. We would like to mention three of them here since they seem
rather glaringly open.

\begin{itemize}

\item Are the results on quantum completely integrable systems
valid under the weaker assumption of classical complete
integrability?  E.g. if $\Delta$ is a quantum integrable
Laplacian, and if its eigenfunctions have a given growth rate, can
the additition of a lower order term such as $\Delta + V$ with $V
\in C^{\infty}(M)$ change the growth rate?

\item Is the existence of a stable elliptic closed geodesic
sufficient for the existence of a sequence $\{\phi_{\lambda}\}$ of
eigenfunctions satisfying $||\phi_{\lambda}||_{\infty} \sim
\lambda^{\frac{n - 1}{8}}?$ Note that if would be sufficient to
obtain sup norm estimates for KAM systems.

\item How do the results generalize to Schroedinger operators
$\Delta + V$ on $\R^n?$ Here, we may assume $V$ grows at infinity
to ensure that the spectrum is discrete.

\end{itemize}

\end{document}